\documentclass[a4paper,12p]{article}
\usepackage{amsmath,amsfonts,amsthm,amssymb,amscd,latexsym, bm,bbm,cite}
\usepackage{hyperref}
\usepackage[usenames,dvipsnames]{color}
\usepackage{authblk}
\hypersetup{colorlinks,
	citecolor=Blue,
	linkcolor=Blue,
	urlcolor=Blue}

\usepackage[T2A]{fontenc}
\usepackage{graphicx}
\usepackage{mathrsfs}

\usepackage{enumitem}

{\theoremstyle {definition} \newtheorem {defi} {Definition} [section] }
{\theoremstyle {plain}  \newtheorem {theo} [defi] {Theorem}}
{\theoremstyle {plain}  \newtheorem {cor} [defi]{Corollary}}
{\theoremstyle {plain}  \newtheorem {lem} [defi]{Lemma}}
{\theoremstyle {plain}  \newtheorem {prop} [defi]{Proposition}}
{\theoremstyle {plain} \newtheorem {rem}[defi] {Remark}}
{\theoremstyle {remark} }

\numberwithin{equation}{section}
\newtheorem*{defi*}{Definition}
\newtheorem*{problem*}{Problem}

\newtheorem*{rem*}{Remark}
\newtheorem*{note*}{Note}

\makeatletter \@addtoreset{equation}{section} \makeatother

\newcommand{\mR}{\mathbb{R}}

\newcommand{\mN}{\mathbb{N}}

\newcommand{\mI}{\mathbb{I}}
\newcommand{\mP}{\mathbb{P}}

\newcommand{\BB}{{\cal B}}
\newcommand{\CC}{{\cal C}}
\newcommand{\DD}{{\cal D}}

\newcommand{\HH}{{\cal H}}

\newcommand{\PP}{{\cal P}}

\newcommand{\VV}{{\cal V}}

\newcommand{\RR}{{\cal R}}

\newcommand{\PPP}{{\cal P}}

\newcommand{\eps}{\varepsilon}
\newcommand{\ph}{\varphi}

\newcommand{\al}{\alpha}

\newcommand{\de}{\delta}
\newcommand{\De}{\Delta}

\newcommand{\supp}{\operatorname{supp}}

\newcommand{\Var}{\operatorname{Var}}

\newcommand{\ov}{\overline}
\newcommand{\wid}{\widetilde}

\newcommand{\ssk}{\smallskip}

\newcommand{\bsk}{\bigskip}
\newcommand{\chp}{\partial}
\newcommand{\MO}{\mathbb{E}}

\newcommand{\ds}{\displaystyle{}}
\newcommand{\ra}{\rightarrow}

\newcommand{\fr}{\frac}

\newcommand{\qmb}{\quad\mbox}
\newcommand{\qu}{\quad}
\newcommand{\qnd}{\quad\mbox{and}\quad}

\newcommand{\sli}{\sum\limits}
\newcommand{\ili}{\int\limits}
\newcommand{\ilif}{\ili_{-\infty}^\infty}
\newcommand{\lbl}{\label}

\newcommand{\llim}{\lim\limits}

\newcommand{\ass}{\quad\mbox{as}\quad}

\newcommand{\rprop}{Proposition \nolinebreak}
\newcommand{\rtheo}{Theorem \nolinebreak}
\newcommand{\rlem}{Lemma \nolinebreak}
\newcommand{\rcor}{Corollary \nolinebreak}

\newcommand{\rsec}{Section \nolinebreak}

\newcommand{\bee}{\begin{equation}}
\newcommand{\eee}{\end{equation}}
\newcommand{\btt}{\begin{theo}}
	\newcommand{\ett}{\end{theo}}
\newcommand{\bl}{\begin{lem}}
	\newcommand{\el}{\end{lem}}
\newcommand{\bpp}{\begin{prop}}
	\newcommand{\epp}{\end{prop}}
\newcommand{\bcc}{\begin{cor}}
	\newcommand{\ecc}{\end{cor}}
\newcommand{\bdd}{\begin{def}}
	\newcommand{\edd}{\end{def}}
\newcommand{\brr}{\begin{rem}}
	\newcommand{\err}{\end{rem}}

\newcommand{\non}{\nonumber}
\newcommand{\sck}{\substack}

\newcommand{\Conf}{\mbox{Conf}}
\newcommand{\CoR}{\mbox{Conf}(\mR)}
\newcommand{\CoRd}{\mbox{Conf}(\mR^d)}
\newcommand{\sCoR}{\small{\mbox{Conf}(\mR)}}
\newcommand{\sCoRd}{\small{\mbox{Conf}(\mR^d)}}
\newcommand{\CCP}{\CC_\mP}

\newcommand\Camp{\mathcal{C}_{\mathbb{P}}}
\newcommand\logd{d^{\mathbb{P}}}

\title{The logarithmic derivative for point processes with equivalent Palm measures}
\author[1,2,3,4,5]{Alexander I. Bufetov\thanks{bufetov@mi.ras.ru}}
\author[2]{Andrey V. Dymov\thanks{dymov@mi.ras.ru}}
\author[6]{Hirofumi Osada\thanks{osada@math.kyushu-u.ac.jp}}
\affil[1]{Aix-Marseille Universit\'e, CNRS, Centrale Marseille, I2M, UMR 7373, 39 rue F. Joliot Curie, Marseille, France}
\affil[2]{Steklov Mathematical Institute of RAS, Moscow, Russia}
\affil[3]{National Research University Higher School of Economics, Moscow, Russia}
\affil[4]{Institute for Information Transmission Problems, Moscow, Russia}
\affil[5]{The Chebyshev Laboratory, Saint-Petersburg State University, Saint-Petersburg, Russia}
\affil[6]{Kyushu University, Faculty of Mathematics, Fukuoka, Japan}
\date{}

\begin{document}
	
	\large\maketitle
	
	\begin{abstract}
		The {\it logarithmic derivative} of a point process plays a key 
		r{\^o}le
		in the general approach, due to the third author,  to constructing  
		diffusions preserving a given point process.
		In this paper we explicitly compute the logarithmic derivative for determinantal 
		processes on $\mR$ with integrable kernels, a large class that includes 
		all the classical processes of random matrix theory as well as processes 
		associated with de Branges spaces.
		The argument uses the quasi-invariance of our processes established by 
		the first author.
	\end{abstract}
	
	\section{Introduction}
	
	Let $\mP$ be a point process on $\mR^d$, or, in other words,  a Borel probability measure on the space of locally finite configurations $\Conf(\mR^d)$.
	It is a natural question whether  one can construct a diffusion $\xi(t)=(\xi^1(t),\xi^2(t),\ldots, \xi^i(t), \dots )$ on the space $\big(\mR^d\big)^\mN$ such that
	the configuration 
	$X(t)=\{\xi^1(t),\xi^2(t),\ldots, \xi^i(t), \dots\}$ 
	is almost surely locally finite for every $t\in {\mathbb R}_+$, and the process $X(t)$, considered as a process on the space $\Conf(\mR^d)$, preserves
	the measure $\mP$.
	For example, if
	$ \mathbb{P} $  is the standard Poisson point process on $\mR^d$,
	then $\xi^i(t)$ are  independent Brownian motions.
	In the series of papers \cite{k-o.fpa, o-h.bes, o.isde,	o.rm, o.rm2, o-t.sm, o-t.tail, o-t.airy} the third author with collaborators developed a general  approach  to constructing the process $\xi$.
	The key step is the computation of the
	{\it logarithmic derivative} $ \logd $ of the measure $ \mathbb{P} $,
	$\logd:\,\mR^d\times\CoRd\mapsto \mR^d$, introduced by the third author in \cite{o.isde}.
	The process $\xi$ is then  a solution of the infinite-dimensional stochastic differential equation
	\bee\non
	\xi ^i (t) = \xi^i (0) +
	B^i(t) +
	\frac{1}{2} \int_0^t
	\logd (\xi ^i(u), X_i(\xi(u)))du,
	\qu i\in\mN,
	\eee
	where the configuration $X_i$ is defined by the formula
	$X_i(\xi(u)):=\{\xi^j(u)\}_{j\neq i}$ and $B^i$ are independent Brownian motions.
	In \cite{o.isde,o-t.airy,o-h.bes}
	logarithmic derivatives were calculated for determinantal processes arising in random matrix theory:
	sine$_2 $, Airy$_2 $, Bessel$_2 $ and the Ginibre point processes.
	The computation was based on  finite particle approximation and had to be
	adapted for each determinantal process separately.
	
	Theorem \ref{thm:logder},  the main result of this paper,  establishes existence  and gives an explicit formula for the logarithmic derivative for  determinantal point processes on $ \mathbb{R}$ with integrable kernels studied in \cite{Buf}, a class that includes, in particular, determinantal processes mentioned above and those corresponding to de Branges spaces \cite{buf-shirai} .
	
	There are other methods to constructing infinite-dimensional diffusions.
	In particular, in \cite{KT07b,KT11}, using extended determinantal kernels, 
	Katori and Tanemura constructed diffusions reversible with respect to the sine$ _2$, Airy$_2 $,
	and Bessel$_2 $ point processes. 
	A different approach to studying the diffusion preserving the sine$_2$ process is due to
	L.-C. Tsai \cite{tsai}. In \cite{BO-diff}, Borodin and Olshanski gave a construction of infinite-dimensional diffusions as scaling limits of random walks on partitions.

	To explain our results in more details we first give an informal definition of the logarithmic derivative. Consider a point process $\mP$ on $\mR^d$
	which admits a differentiable first correlation function $\rho_1:\mR^d\mapsto\mR$.
	Denote by $ \mathbb{P}^a $ the reduced Palm measure conditioned at the point $ a \in \mathbb{R}^d $ and define the reduced Campbell measure $ \Camp $ as a Borel sigma-finite measure on the space
	$\mR^d\times \mathrm{Conf} (\mathbb{R}^d)$
	given by
	\bee\non
	d\Camp (a,X) = \rho_1 (a) d\mathbb{P}^a (X)da.
	\eee
	Then, informally, the logarithmic derivative $d_\mP$
	is defined
	as a gradient of the logarithm of $ \Camp $,
	\begin{align}\label{:15}&
	\logd (a,X) = \nabla_a \big(\ln \rho_1(a) +
	\ln \mathbb{P}^a (X)\big)
	,\end{align}
	see Definition \ref{def:logder}.
	The main problem when proving the existence of the logarithmic derivative is to give a sense to the term
	$ \nabla_a \ln \mathbb{P}^a (X)$.
	
	Our first result is Proposition \ref{lem:logder},
	where we find the connection between equivalence of the Palm measures conditioned at different points and the existence of the logarithmic derivative.
	More specifically, we consider a point process $\mP$ on $\mR^d$ as above and assume that for any $a,b\in\mR^d$ the reduced Palm measures
	$\mP^a$ and $\mP^b$ are equivalent,
	$$
	d\mP^b(X)=\RR_{b,a}(X)d\mP^a(X),
	$$
	the Radon-Nikodym derivative $\RR_{b,a}(X)$ is
	continuous with respect to $b$ in $L^1(\mP^a,\CoRd)$ and
	the derivative $\nabla_b\RR_{b,a}$ exists
	in appropriate sense.
	We then prove that
	the logarithmic derivative $\logd$ exists and the formula \eqref{:15} is valid
	with
	\bee\non
	\nabla_a\ln \mathbb{P}^a:= \nabla_b\big|_{b=a}\RR_{b,a}.
	\eee
	Our second and main result is the mentioned above \rtheo\ref{thm:logder}.
	To establish it, it suffices to check that
	assumptions of \rprop\ref{lem:logder}
	are satisfied for the
	considered class of determinantal point processes.
	To show this, we use the results of the paper \cite{Buf},
	where the first author proved that for this class of determinantal processes the reduced Palm measures are equivalent and the Radon-Nikodym derivative has the form
	$$
	\RR_{b,a}=\llim_{R\ra\infty, \de\ra 0 } \RR_{b,a}^{R,\de}
	\qmb{where}\quad
	\RR_{b,a}^{R,\de} = C^{R,\de}_{b,a}
	\prod\limits_{\substack{X\in\sCoR:\, |x|<R,\\ |x-a|,|x-b|>\de}}\left(\frac{x-b}{x-a}\right)^2,
	$$
	and $C^{R,\de}_{b,a}$ are some normalizing constants.
	While the continuity in $b$ of $\RR_{b,a}$ follows immediately from the results of
	\cite{Buf}, the proof of its differentiability
	requires some efforts.
	To get it, we approximate the Radon-Nikodym derivative
	$\RR_{b,a}$ by the function $\RR_{b,a}^{R,\de}$,
	compute the derivative of the latter,
	and then pass to the limit $R\ra\infty$, $\de\ra 0$
	using the techniques of normalized additive and multiplicative
	functionals developed in \cite{Buf},
	which we outline in the appendix.
	Finally, we find
	$$
	\nabla_b\big|_{b=a}\RR_{b,a}
	=\llim_{R\ra\infty, \de\ra 0} (S_a^{R,\de}-\MO^a S_a^{R,\de}),
	$$
	where
	$\ds{S_a^{R,\de}=\sli_{\substack{x\in X:\, |x|<R, \\ |x-a|>\de}}}\fr{2}{a-x}$
	and $\MO^a$ stands for the expectation with respect to the reduced Palm measure $\mP^a.$

	The paper is organized as follows.
	In Section 2 we  formulate our main results,
	\rprop \ref{lem:logder} and \rtheo \nolinebreak\ref{thm:logder}.
	Section 3 is devoted to the proofs.
	In Appendix \nolinebreak\ref{sec:regularization}
	we recall some results of \cite{Buf}
	needed in the proof of \rtheo \nolinebreak\ref{thm:logder}.

	\section{Formulation of the main results}
	
	\subsection{Configurations, point processes, Palm distributions}
	\lbl{sec:setup}
	
	Consider the space of locally finite configurations
	\bee\non
	\CoRd:=\big\{X\subset\mR^d |\, X \mbox{ does not have limit points in }\mR^d \big\}.
	\eee
	A Borel probability measure $\mP$ on $\CoRd$ is called a point process.
	Take a bounded Borel set $B\in\BB(\mR^d)$ and consider a function
	$\#_B:\Conf(\mR^d)\mapsto\mN\cup\{0\}$ such that $\#_B(X)$ is equal to the cardinality of the set $B\cap X$.
	Assume that the process $\mP$ admits the first correlation function $\rho_1$, that is
	for any bounded $B\in\BB(\mR^d)$
	the function $\#_B$ is integrable with respect
	to the measure $\mP$ and there exists
	a function $\rho_1\in L^1_{loc}(\mR^d)$ satisfying
	$$
	\int_{B}\rho_1(x)\, dx= \int_{\small\Conf(\mR^d)} \#_B(X)d\mP(X), \quad \forall B\in\BB(\mR^d),
	\qu B\mbox{ is bounded}
	.
	$$
	Define the first correlation measure $\hat\rho_1$ as $\ds{\hat\rho_1(B)=\int_{B}\rho_1(x)\, dx}$.
	
	The {\it Campbell measure}  $\hat\CCP$
	is a sigma-finite Borel measure on $\mR^d\times \Conf(\mR^d)$ defined as
	$$
	\hat\CCP(B, {\mathscr Z})=\displaystyle \int\limits_{{\mathscr Z}} \#_B(X)d\mP(X),
	\qu \forall B\in\BB(\mR^d), \qu {\mathscr Z}\in\BB\big(\Conf(\mR^d)\big),
	$$
	where $\BB\big(\Conf(\mR^d)\big)$ stands for the Borel sigma-algebra on $\Conf(\mR^d)$.
	Fix a Borel set ${\mathscr Z}\subset \Conf(\mR^d)$
	and consider a sigma-finite measure $\CCP^{{\mathscr Z}}$
	on $\mR^d$ given by the formula
	$$
	\CCP^{{\mathscr Z}}(B)=\hat\CCP(B,{\mathscr Z}), \qu \forall B\in\BB(\mR^d).
	$$
	By definition, for any ${\mathscr Z}\in \BB\big(\Conf(\mR^d)\big)$ the measure
	$\CCP^{{\mathscr Z}}$ is absolutely continuous
	with respect to $\hat\rho_1$.
	Then the {\it Palm measure} ${\hat \mP}^a$,
	defined for $\hat\rho_1$-almost every $a\in \mR^d$,
	is a measure on $\Conf(\mR^d)$
	given by the relation
	$$
	{\hat \mP}^a({\mathscr Z})=\displaystyle \frac{d\CCP^{{\mathscr Z}}}{d\hat\rho_1}\left(a\right).
	$$
	Equivalently, the Palm measure ${\hat \mP}^a$ is the canonical conditional measure of the Campbell measure $\hat\CCP$
	with respect to the measurable partition of the space $\mR^d\times \Conf(\mR^d)$
	into subsets of the form $\{a\}\times \Conf(\mR^d)$, $a\in \mR^d$.
	Thus, we can write
	$$
	d\hat\CCP(a,X)=d\hat\mP^a(X) \rho_1(a)\,da.
	$$
	By definition, the Palm measure ${\hat \mP}^{a}$ is supported on the subset of configurations containing a particle at the position $a$.
	The {\it reduced
		Palm measure} $\mP^{a}$
	is defined
	as the push-forward of the Palm measure
	${\hat \mP}^{a}$ under the map
	$X\to X\setminus \{a\}$ erasing the particle $a$ from the configuration $X$.
	We then define the {\it reduced Campbell measure} $\CCP$ as
	
	\bee\lbl{Campblm}
	d\CCP(a,X)=d\mP^a(X) \rho_1(a)\, da.
	\eee
	Note that, in difference with the notions of the (reduced) Palm measure and the Campbell measure,
	this definition is not standard. Writing it in a more formal way, we obtain
	$$
	\CCP(B,\mathscr{Z})=\int_B\int_\mathscr{Z} d\mP^a(X) \rho_1(a)\, da
	\qu \forall B\in\BB(\mR^d), \qu {\mathscr Z}\in\BB\big(\Conf(\mR^d)\big).
	$$
	For more details see e.g. \cite{Buf} and \cite{DVJ}.

	\subsection{Definition of the logarithmic derivative} A function $\ph:\Conf(\mR^d)\mapsto\mR$ is called {\it local}
	if there exists a  compact set $K\subset \mR^d$ such that
	$\ph(X)\equiv\ph(X\cap K)$.
	For a local function $\ph$
	we define symmetric functions
	$\ph_n: \mR^{nd}\mapsto \mR$, $n\geq 1$,
	by the relation
	\bee\non
	\ph_n(x_1,\ldots,x_n)=\ph\big(\{x_1,\ldots,x_n\}\big).
	\eee
	We say that a local function $\ph$ is smooth if
	the functions $\ph_n$ are smooth for all $n\in\mN$.
	We denote by $\DD_0$ the space of all bounded local smooth functions on $\CoRd$.
	
	Denote by $B_R$ a ball in $\mR^d$ of radius $R$.
	Let
	$L_{loc}^1(\mR^d\times\CoRd,\CCP)$
	be the space of vector-functions
	$f:\mR^d\times\CoRd \mapsto \mR^d$
	satisfying
	$f\in L^1(B_R\times\CoRd,\CCP)$, for all $R>0$.
	
	Denote by $C_0^\infty$ the space of smooth real-valued functions on $\mR^d$ which have compact supports. 
	We say that a function $\ph:\mR^d\times\CoRd\mapsto\mR$ belongs to the space $C_0^\infty \DD_0$ if $\ph=\ph_1\ph_2$, where $\ph_1\in C_0^\infty(\mR^d)$ and $\ph_2\in\DD_0$.
	\begin{defi}
		\lbl{def:logder}
		Let $\mP$ be  a point process  on $\mR^d$ that admits the first correlation function.
		A function $d^\mP\in L_{loc}^1(\mR^d\times\CoRd,\CCP)$ is called the {\it logarithmic derivative} of $\mP$ if
		for any observable $\ph\in C_0^\infty\DD_0$  we have
		\bee\non
		\ili_{\mR^d\times\small{\CoRd}} \nabla_a\ph(a,X) \,d\CCP(a,X)
		=-\ili_{\mR^d\times{\small \CoRd}} d^\mP(a,X) \ph(a,X) \,d\CCP(a,X).
		\eee
	\end{defi}
	
	\subsection{Logarithmic derivative of a point process with equivalent Palm measures}
	
	Consider a point process $\mP$ on $\mR^d$ that admits the first correlation function $\rho_1$; recall that we denote by $\hat{\rho}_1$ the first correlation measure, 
	$$\hat{\rho}_1(da)=\rho_1(a)\,da.$$
	In this subsection we give a general scheme for the computation of the logarithmic derivative $\logd$ under the following assumption.

	\smallskip
	
	{\bf Assumption 1.}
	\begin{enumerate}
		\item
		The first correlation function $\rho_1$ is $C^1$-smooth.
		\item
		For $\hat{\rho}_1$-almost all $a,b\in\mR^d$ the reduced Palm measures $\mP^a$ and $\mP^b$ are equivalent.
	\end{enumerate}
	Denote by $\RR_{b,a}$ their Radon-Nikodym derivative, so that
	$$
	d\mP^b(X)=\RR_{b,a}(X)\,d\mP^a(X).
	$$
	\begin{enumerate}[resume]
		\item
		For $\hat{\rho}_1$-almost all $a\in\mR^d$ we  have $\RR_{b,a}\ra 1$ as $b\ra a$ in
		$L^1(\mP^a,\CoRd)$.
	\end{enumerate}
	For a function $\ph\in C^\infty_0(\mR^d)\DD_0$ we define the function $f_\ph:\, \mR^d\mapsto\mR$ as
	$$
	f_\ph(\eps):=\ili_{\mR^d\times{\small \CoRd}} \RR_{a+\eps,a}(X)\ph(a,X)\,d\CCP(a,X).
	$$
	\begin{enumerate}[resume]
		\item
		For any $\ph\in C^\infty_0(\mR^d)\DD_0$ the function
		$f_\ph$ admits partial derivatives in $\eps$ at the point $\eps=0$.
		There exist functions $\chp_{i} \RR:\mR^d\times \CoRd\mapsto\mR$
		such that for any $\ph$ as above and any $1\leq i\leq d$ we have
		$$
		\chp_{\eps_i}f_{\ph}(0)
		=\ili_{\mR^d\times{\small \CoRd}}
		\chp_{i}\RR(a,X)\ph(a,X)\,d\CCP(a,X).
		$$
	\end{enumerate}
	Set
	$$
	\nabla\RR:=(\chp_{1}\RR, \ldots, \chp_{d}\RR ).
	$$
	\bpp\lbl{lem:logder}
	Let $\mP$ be a point process on $\mR^d$ satisfying Assumption 1.
	Then for $\hat{\rho}_1$-almost all $a\in\mR$ its logarithmic derivative $d_{\mP}$ exists and has the form
	\bee\lbl{logder_gener}
	d_{\mP}(a,X)=\nabla_a\ln\rho_1(a) + \nabla\RR(a,X).
	\eee
	\epp
	Note that there is no need to define the logarithmic derivative at the points $a\in\mR^d$ where $\rho_1(a)=0$ since the measure $\hat{\rho}_1$ of the set of such $a$ is zero. 
	Proof of \rprop \ref{lem:logder} is given in Section \ref{sec:lemlogderpr}. 
	It is based on the differentiation by parts, that is why we crucially need  the absolute continuity of the measure $\hat{\rho}_1$ and the differentiability of its density, 
	which is the first correlation function $\rho_1$.

	\subsection{Logarithmic derivative of a determinantal process on $\mathbb{R}$ with an integrable kernel}
	
	In this section we construct the logarithmic derivative for a
	class of determinantal processes on $\mR$.
	A point process
	$\mP$ on $\CoR$
	is called {\it determinantal} if there exists a locally
	trace class operator
	$\PP: L^2(\mR,dx)\mapsto L^2(\mR,dx)$
	such that for any bounded measurable function $h$,
	for which the support
	$\supp (h-1)=:D$ is a compact set,
	we have
	$$
	\MO\left(\prod\limits_{x\in X}h(x)\right)
	=\det\big(1+(h-1)\PP\mI_{D} \big).
	$$
	Here the expectation is taken with respect to the measure $\mP$,
	$\det$ stands for the Fredholm determinant
	and $\mI_D$ denotes the indicator function of the set $D$.
	See for details \cite{Sos00,ST}.
	Since the operator
	$\PP$ is locally trace class,
	it admits a kernel which we denote by $\Pi$.
	We impose the following restrictions for $\PP$ and $\Pi$.
	
	\ssk
	{\bf Assumption 2}.
	
	\begin{enumerate}
		\item
		The operator $\PP$ is an orthogonal projection onto a closed subspace  $L\subset L^2(\mR,dx)$.
		\item
		For $\hat{\rho}_1$-almost all $a\in\mR$,
		given any function $\ph\in L$ satisfying 
		$\ph(a)=0$, we have $(x-a)^{-1}\ph\in L$.
		\item
		The kernel $\Pi$ is $C^2$-smooth on $\mR^2$.
		\item
		We have
		$\ds{\ili_\mR \fr{\Pi(x,x)}{1+x^2} \, dx<\infty}$.
	\end{enumerate}
	Also, note that for any $a\in\mR$ the function $\Pi(a,\cdot)$ belongs to $L^2(\mR,dx)$.
	
	\ssk
	Take $a\in\mR$, $R\gg 1$ and $\de\ll 1$, and consider  the additive functional
	\bee\lbl{S}
	S^{R,\de}_{a}:\Conf(\mR)\mapsto\mR,
	\qu
	S^{R,\de}_{a} (X)
	=\sli_{\sck{x\in X:\, |x|<R, \\ |x-a|>\de}}  \fr{2}{a-x}.
	\eee
	The additive functional $S^{R,\de}_a$ may diverge
	as $R\ra\infty$.
	To overcome this difficulty we define the normalized
	additive functional
	$$
	\ov S^{R,\de}_{a} := S^{R,\de}_{a}-\MO^a S^{R,\de}_{a},
	$$
	where $\MO^a$ stands for the expectation
	with respect to the reduced Palm measure $\mP^a$.
	Results obtained in \cite{Buf} imply that,
	under Assumption 2,
	for $\hat{\rho}_1$-almost all $a\in\mR$ there exists a function
	$\ov S_{a}:\CoR\mapsto\mR$,
	such that
	\bee\lbl{S1conv}
	\ov S^{R,\de}_{a} \ra \ov S_{a} \qmb{as}\qu R\ra\infty,\, \de\ra 0 \qmb{in}\qu L^2(\CoR,\mP^a).
	\eee
	Moreover,
	the convergence (\ref{S1conv}) holds uniformly in
	$a$ as $a\in\mR$ ranges in a compact set.
	The required theory from \cite{Buf} is recalled in Appendix \ref{sec:add}
	and the convergence \eqref{S1conv} is established in
	\rcor\ref{cor:add}.

	\btt\lbl{thm:logder}
	Let $\mP$ be a determinantal process on $\mR$ satisfying Assumption 2.
	Then for $\hat{\rho}_1$-almost all $a\in\mR$ the logarithmic derivative $d_{\mP}$ exists and has the form
	\bee\non
	d_{\mP}(a,X)=\fr{d}{da}\ln\rho_1(a) + \ov S_a(X).
	\eee
	\ett
	\rtheo\ref{thm:logder} is proven in \rsec\ref{sec:thmlogredpr}.
	There, using results of \cite{Buf},
	we show that Assumption 2 implies Assumption 1
	with $\nabla R= \ov S_a$.
	Then  \rtheo\ref{thm:logder} follows from  \rprop\ref{lem:logder}.
	\section{Proofs of the main results}
	
	\subsection{Proof of Proposition \ref{lem:logder}}
	\lbl{sec:lemlogderpr}
	
	Take a function $\ph\in C_0^\infty\DD_0$. We have
	$$
	I:=-\ili_{\mR^d\times \sCoRd} \chp_{a_i}\ph(a,X) \,d\CCP(a,X)
	= -\ili_{\mR^d\times \sCoRd} \lim\limits_{\eps\ra 0}\fr{\ph(a+\eps_i ,X)-\ph(a,X)}{\eps} \,d\CCP(a,X),
	$$
	where $\eps_i:=\eps e_i$ and $e_i$ is the $i$-th basis vector of $\mR^d$.
	Using the dominated convergence theorem,
	we exchange the limit with the integral.
	The latter can be applied since
	$$
	\Big|
	\fr{\ph(a+\eps_i,X)-\ph(a,X)}{\eps}
	\Big|
	\leq \sup\limits_{X\in\sCoRd,\, x\in\mR^d}
	\big|\chp_{x_i}\ph(x,X)\big|
	$$
	and $\ph\in C_0^\infty\DD_0$. We get
	\begin{align}\non
	I&=-\llim_{\eps\ra 0}\fr{1}{\eps}
	\Big(
	\ili_{\mR^d\times \sCoRd} \ph(a+\eps_i,X) \,d\CCP(a,X)
	-\ili_{\mR^d\times \sCoRd} \ph(a,X) \,d\CCP(a,X)
	\Big)
	\\ \lbl{I1}
	&=\llim_{\eps\ra 0}\fr{1}{\eps}
	\Big(\ili_{\mR^d\times \sCoRd} \ph(a,X) \,d\CCP(a+\eps_i,X)
	-\ili_{\mR^d\times \sCoRd} \ph(a,X) \,d\CCP(a,X)
	\Big),
	\end{align}
	where in the last line of (\ref{I1})
	we put $\eps:= -\eps$.
	Using the definition of the reduced Campbell measure (\ref{Campblm}),
	we find
	\begin{align}\non
	I&=\llim_{\eps\ra 0}\fr{1}{\eps}\Big(\ili_{\mR^d\times \sCoRd}
	\ph(a,X) \rho_1(a+\eps_i)\,d\mP^{a+\eps_i}(X)da
	\\\non
	&- \ili_{\mR^d\times \sCoRd}
	\ph(a,X) \rho_1(a)\,d\mP^a(X)da\Big)
	\\\non
	&=\llim_{\eps\ra 0}\fr{1}{\eps}\Big(\ili_{\mR^d\times \sCoRd}
	\ph(a,X) \big( \rho_1(a+\eps_i)-\rho_1(a)\big)\,d\mP^{a+\eps_i}(X)da
	\\\non
	&+ \ili_{\mR^d\times \sCoRd}
	\ph(a,X) \rho_1(a)\,\big(d\mP^{a+\eps_i}(X)-d\mP^a(X)\big)da\Big)
	\\\nonumber
	&=\llim_{\eps\to 0}(I_1^\eps+I_2^\eps).
	\end{align}
	Using Assumption 1(3), we obtain
	\bee\lbl{I_1^eps}
	\llim_{\eps\to 0} I_1^\eps
	= \ili_{\mR^d\times \sCoRd} \ph(a,X) \chp_{a_i}\rho_1(a) \,d\mP^a(X)da
	=\ili_{\mR^d\times \sCoRd} \ph(a,X) \chp_{a_i}\big(\ln\rho_1(a)\big) \,d\CCP(a,X).
	\eee
	In view of Assumption 1(2), we have
	\bee\non
	\llim_{\eps\to 0} I_2^\eps
	=\llim_{\eps\to 0}\fr{1}{\eps}\Big(
	\ili_{\mR^d\times \sCoRd} \ph(a,X) \RR_{a+\eps_i,a}(X) \,d\CCP(a,X)
	-\ili_{\mR^d\times \sCoRd} \ph(a,X)\,d\CCP(a,X) \Big).
	\eee
	Then, because of the identity $\RR_{a,a}(X)\equiv 1$,
	Assumption 1(4) implies
	\bee\lbl{I_2^eps}
	\llim_{\eps\to 0} I_2^\eps
	=\ili_{\mR^d\times \sCoRd} \ph(a,X) \chp_i \RR(a,X) \,d\CCP(a,X).
	\eee
	Combining \eqref{I_1^eps} with \eqref{I_2^eps}
	we obtain the desired relation \eqref{logder_gener}.
	\qed

	\subsection{Proof of Theorem \ref{thm:logder}}
	\lbl{sec:thmlogredpr}
	
	In view of  \rprop \ref{lem:logder},
	it suffices to check that Assumption 1 is satisfied with
	$
	\nabla\RR(a,X)=\ov S_a(X).
	$
	Item 1 of Assumption 1 immediately follows from item 2 of Assumption \nolinebreak 2.
	The proof of the other items relies on the results obtained in the paper \cite{Buf};
	see also \cite{Buf1}.
	One of the main tools we use borrowed from the works above is the following lemma.
	Take $a,b\in\mR$ and
	consider the normalized multiplicative functional
	$\ov\Psi_{b,a}^{R,\de}: \CoR\mapsto \mR$ given by
	\bee\lbl{Psi}
	\ov\Psi_{b,a}^{R,\de}(X):= C^{R,\de}_{b,a}
	\prod_{\sck{x\in X: \, |x|<R, \\ |x-a|, |x-b|>\de}}
	\Big(
	\fr{x-b}{x-a}
	\Big)^2,
	\eee
	where the constant $C^{R,\de}_{a,b}$
	is specified by the normalization requirement
	$\MO^a \ov\Psi_{b,a}^{R,\de} = 1$.
	Here and further on we set
	$\prod\limits_{x\in\emptyset} f(x)=1$, for any function $f$.
	\bl\lbl{lem:multf}
	{\bf 1.}
	Under Assumption 1, there exists $\al>0$ and a function
	$\ov\Psi_{b,a}: \mathrm{Conf}(\mR)\mapsto \mR$
	satisfying
	\bee\non
	\ov\Psi_{b,a}^{R,\de}\ra \ov\Psi_{b,a} \qmb{as}\qu R\ra\infty,\,\de\ra 0
	\qmb{in}\qu  L^{1+\al}(\mathrm{Conf}(\mR),\mP^a),
	\eee
	for $\hat{\rho}_1$-almost all $a\in\mR$, uniformly in $a,b$ which range in compact subsets of $\mR$.
	
	{\bf 2.}
	For $\hat{\rho}_1$-almost all $a,b\in\mR^d$, the function $\ov\Psi_{b,a}$
	is the Radon-Nikodym derivative of the Palm measure $\mP^b$
	with respect to the the Palm measure $\mP^a$, i.e.
	\bee\non
	d\mP^b (X) = \ov\Psi_{b,a}(X) \,d\mP^a(X).
	\eee
	\el
	{\it Proof.}
	Item (1) is established in \rcor\ref{cor:mult}(1) and follows from
	results obtained in \cite{Buf},
	which we explain in Appendix \ref{sec:mult}.
	Item \nolinebreak (2) is a
	particular case
	of
	\rtheo 1.4\,(1) from \cite{Buf}.
	Note that in  \cite{Buf} the multiplicative functional
	is defined
	as the product over the set
	$\{x\in X:\, |x|<R, |x-a|>\de\}$,
	so that in difference
	with the definition (\ref{Psi})
	the point $b$ is not isolated.
	It can be checked directly that
	this does not affect the proof at all.
	\qed
	
	\ssk
	Lemma \ref{lem:multf}(2) implies item 2 of Assumption \nolinebreak 1
	with $\RR_{b,a}=\ov\Psi_{b,a}$.
	Because of the bounds $|x|<R$ and $|x-a|>\de$,
	the functions
	$\ov\Psi_{b,a}^{R,\de}(X)$
	are $\mP^a$-almost surely continuous with respect to $b$.
	Then, using the dominated convergence theorem,
	we see that they are continuous in $L^1(\CoR,\mP^a)$.
	Then the uniformity in $b$ of convergence from Lemma \ref{lem:multf}(1)
	implies that the functions $\ov\Psi_{b,a}$ also are continuous with respect to $b$ in $L^1(\CoR,\mP^a)$.
	So that item 3 of Assumption \nolinebreak 1 is fulfilled as well.
	
	It remains to check that item 4 of Assumption \nolinebreak 1 holds with
	$\nabla \RR(a,X)=\ov S_a(X)$.
	Due to \rlem\ref{lem:multf}(2), we have
	$$
	f_\ph(\eps) = \ili_{\mR\times\sCoR}
	\ov\Psi_{a+\eps,a}(X)\ph(a,X)\,d\CCP(a,X).
	$$
	We need to show that the function $f_\ph$ is differentiable
	at zero and
	\bee\lbl{I'2des}
	\fr{d}{d\eps}f_\ph(0)=\ili_{\mR\times\sCoR}
	\ov S_{a}(X)
	\ph(a,X)\,d\CCP(a,X).
	\eee
	Due to Lemma \ref{lem:multf}(1), we have
	\bee\lbl{mmain}
	f_\ph(\eps)=\llim_{R\ra\infty,\de\ra 0} f_\ph^{R,\de}(\eps),
	\eee
	where
	$$
	f_\ph^{R,\de}(\eps) = \ili_{\mR\times\sCoR}
	\ov\Psi_{a+\eps,a}^{R,\de}(X)\ph(a,X)\,d\CCP(a,X).
	$$
	\bpp\lbl{lem:chpeps}
	The function $f_\ph^{R,\de}$
	is $C^1$-smooth.
	For any  $\eps$ from a small neighbourhood of zero
	its derivative
	$\fr{d}{d\eps}f_\ph^{R,\de}(\eps)$ converges as $R\ra\infty, \de\ra 0$,
	uniformly in $\eps\ll 1$.
	\epp
	Proof of \rprop \ref{lem:chpeps} is given in the next subsection.
	Jointly with \eqref{mmain}, \rprop \ref{lem:chpeps} implies that the function $f_\ph$ is differentiable in a small neighbourhood of zero and
	\bee\lbl{formal}
	\fr{d}{d\eps} f_\ph(\eps)
	=\llim_{R\ra\infty,\de\ra 0} \fr{d}{d\eps}f_\ph^{R,\de}(\eps).
	\eee
	We have
	\bee\non
	\fr{d}{d\eps}f_\ph^{R,\de}(\eps)
	=\ili_{\mR\times\sCoR}
	\fr{d}{d\eps}\ov\Psi_{a+\eps,a}^{R,\de}(X)\ph(a,X)\,d\CCP(a,X).
	\eee
	By definition (\ref{Psi}) of the multiplicative functional $\ov\Psi^{R,\de}_{a+\eps,a}$, we have
	\bee\lbl{Psi'}
	\fr{d}{d\eps}\ov\Psi^{R,\de}_{a+\eps,a}
	=\fr{d}{d\eps}
	\exp\Big(
	\ln C^{R,\de}_{a+\eps,a}
	+ 2\sli_{\sck{x\in X:,\, |x|<R, \\ |x-a|,|x-(a+\eps)|>\de}}
	\ln\Big|\fr{x-(a+\eps)}{x-a}\Big|
	\Big)
	=\ov\Psi^{R,\de}_{a+\eps,a} \Big(
	\fr{d}{d\eps}\ln  C^{R,\de}_{a+\eps,a} + S^{R,\de}_{a,a+\eps}
	\Big),
	\eee
	where
	$$
	S^{R,\de}_{a,b}:
	=\sli_{\sck{x\in X:\, |x|<R, \\ |x-a|,|x-b|>\de}}  \fr{2}{b-x}.
	$$
	Since, by definition, $\MO^a \ov\Psi^{R,\de}_{a+\eps,a}\equiv 1$,
	we have
	$$
	\ds{\MO^a \fr{d}{d\eps}\ov\Psi^{R,\de}_{a+\eps,a}=\fr{d}{d\eps}\MO^a\ov\Psi^{R,\de}_{a+\eps,a}=0}.
	$$
	Then, taking the expectation $\MO^a$ of the both sides of (\ref{Psi'}), we  get
	\bee\lbl{lnC}
	\fr{d}{d\eps}\ln  C^{R,\de}_{a+\eps,a}
	= -\MO^a \big(\ov\Psi^{R,\de}_{a+\eps,a} S^{R,\de}_{a,a+\eps}\big).
	\eee
	Now (\ref{Psi'}) together with (\ref{lnC}) implies
	\bee\lbl{I'1}
	\fr{d}{d\eps}f_\ph^{R,\de}(\eps)
	=\ili_{\mR\times\sCoR}
	\ov\Psi^{R,\de}_{a+\eps,a}
	\big( S^{R,\de}_{a,a+\eps}- \MO^a \ov\Psi^{R,\de}_{a+\eps,a}S^{R,\de}_{a,a+\eps}\big)
	\ph(a,X)\,d\CCP(a,X).
	\eee
	Since $\ov\Psi^{R,\de}_{a,a}=1$ and $S_{a,a}^{R,\de}=S_a^{R,\de}$,
	where the additive functional $S_a^{R,\de}$ is defined in (\ref{S}),
	we obtain
	\bee\lbl{I'2}
	\fr{d}{d\eps}f_\ph^{R,\de}(0)
	=\ili_{\mR\times\sCoR}
	\big( S^{R,\de}_{a}- \MO^a S^{R,\de}_{a}\big)
	\ph(a,X)\,d\CCP(a,X).
	\eee
	Due to \eqref{S1conv}, the function
	$S^{R,\de}_{a}- \MO^a S^{R,\de}_{a}=\ov S^{R,\de}_{a}$
	converges to $\ov S_{a}$ in $L^2(\CoR,\mP^a)$ as $R\ra\infty,\, \de\ra 0$, uniformly in
	$\hat{\rho}_1$-almost all $a\in\cup_{X\in\sCoR} \supp \ph(\cdot,X)$, since the latter set is compact.
	Then the right-hand side of \eqref{I'2} converges to that of \eqref{I'2des}.
	In view of \eqref{formal}, this concludes the proof of the theorem.
	\qed

	\subsection{Proof of Proposition \ref{lem:chpeps}}
	The $C^1$-smoothness of the function
	$f_\ph^{R,\de}$ is obvious since
	its derivative has the form
	(\ref{I'1}).
	So that, we only need to establish
	the uniform in $\eps$ convergence of the derivative
	$\ds{\fr{d}{d\eps}f_\ph^{R,\de}}$
	as $R\ra\infty,\de\ra 0$.
	Clearly, it suffices to show that the function
	\bee\lbl{dJ}
	J^{R,\de}(a,b):
	=\ov \Psi_{b,a}^{R,\de}
	\big( S^{R,\de}_{a,b}- \MO^a \ov\Psi^{R,\de}_{b,a} S^{R,\de}_{a,b}\big)
	\qmb{converges as $R\ra\infty,\,\de\ra 0$,}
	\eee
	in $L^1(\CoR, \mP^a)$ uniformly in $b$ and $\hat{\rho}_1$-almost all $a$,
	where $a$ ranges in a compact set
	and $b$ satisfies $|a-b|< \vartheta$,
	for some fixed $\vartheta\ll 1$.
	All convergences below will be uniform in $a,b$ satisfying these restrictions, 
	and we do not mention it any more.
	Further on we assume $R$ to be sufficiently large
	and $\de$ to be sufficiently small
	where it is needed.
	
	For real Borel functions $f,g$ where $g$ is non-negative
	we define the additive and multiplicative functionals
	$S_f,\ov S_f,\Psi_g,\wid\Psi_g,\ov\Psi_g:\,\CoR\mapsto\mR$ by the formulas
	(\ref{addfunc}), (\ref{normaddfunc}),
	(\ref{multfunc}) and (\ref{normmultfunc}).
	Clearly, they are well-defined if the functions $f$
	and $g$ are bounded and the supports $\supp f$, $\supp(g-1)$ are compact.
	However, the normalized functionals $\ov S_f, \wid\Psi_g$ and $\ov\Psi_g$
	can be defined for larger classes of functions, see
	Appendices \nolinebreak \ref{sec:add} and \ref{sec:mult}.
	
	Let us define
	$$
	h^R_>(x)
	:=\fr{2}{b-x}\mI_{\{|x|<R,\,|x-a|>\de,\,|x-b|\geq 1\}}(x)
	=\fr{2}{b-x}\mI_{\{|x|<R,\,|x-b|\geq 1\}}(x),
	$$
	where we have used that $|a-b|\ll 1$, 
	so that the constraint $|x-a|>\de$ holds automatically. 
	Set also
	$$
	h^\de_<(x)
	:=\fr{2}{b-x}\mI_{\{|x|>R,\,|x-a|>\de,\, \de<|x-b|<1\}}(x)
	=\fr{2}{b-x}\mI_{\{|x-a|>\de,\, \de<|x-b|<1\}}(x).
	$$
	Then we have
	$$
	S_{a,b}^{R,\de}=S_{h^R_>}+S_{h^\de_<}.
	$$
	Recall that
	$\MO^a \ov\Psi_{b,a}^{R,\de}=1$.
	Then, subtracting in the brackets of (\ref{dJ}) the term
	$\MO^a S_{h^R_>}$
	and adding the term
	$\MO^a \ov\Psi_{b,a}^{R,\de}\MO^a S_{h^R_>}$, we obtain
	\bee\lbl{JJJ}
	J^{R,\de}(a,b):
	=\ov \Psi_{b,a}^{R,\de}
	(\ov S_{h^R_>}+S_{h^\de_<})
	- \ov\Psi_{b,a}^{R,\de}\MO^a \big(\ov\Psi_{b,a}^{R,\de} (\ov S_{h^R_>}+S_{h^\de_<})\big).
	\eee
	Now, to establish the convergence (\ref{dJ})
	it suffices to show that the functions
	$\ov\Psi_{b,a}^{R,\de} \ov S_{h^R_>}$
	and
	$\ov\Psi_{b,a}^{R,\de} S_{h^\de_<}$
	converge as $R\ra\infty,\, \de\ra 0$ in $L^1(\CoR,\mP^a)$.
	Indeed, in view of
	\rlem \ref{lem:multf}(1), 
	the function $\ov\Psi_{b,a}^{R,\de}$
	converges in $L^1(\CoR,\mP^a)$ itself,
	so that we will see that the second summand 
	from \eqref{JJJ} converges.
	
	{\it Term $\ov\Psi_{b,a}^{R,\de} \ov S_{h^R_>}$.}
	Let
	$$
	h_>(x):=\fr{2}{b-x}\mI_{\{|x-b|\geq 1\}}(x).
	$$
	Due to \rcor\ref{cor:add}(2),
	the additive functional
	$\ov S_{h_>}$
	is well-defined and we have the convergence
	\bee\lbl{to_S}
	\ov S_{h^R_>}\ra \ov S_{h_>}
	\ass
	R \ra \infty
	\qmb{in}\quad
	L^p(\CoR,\mP^a)
	\eee
	with $p=2$.
	We claim that
	it takes place for any $p>2$ as well.
	This concludes consideration of the term
	$\ov\Psi_{b,a}^{R,\de} \ov S_{h^R_>}$
	since, using the H\"older inequality,
	from (\ref{to_S}) joined with
	Lemma \ref{lem:multf}(1)
	we obtain
	$$
	\ov\Psi_{b,a}^{R,\de} \ov S_{h^R_>}
	\ra
	\ov\Psi_{b,a}\, \ov S_{h_>}
	\ass
	\de\ra 0, \,R \ra \infty
	\qmb{in}\quad
	L^1(\CoR,\mP^a).
	$$
	Denote
	$$
	\De^R:=h^R_>-h_>.
	$$
	Due to the Cauchy-Bunyakovsky-Schwarz inequality, we have
	\bee\non
	\MO^a|\ov S_{h_>^R}- \ov S_{h_>}|^p
	= \MO^a|\ov S_{\De^R}|^p
	\leq \sqrt{\MO^a(\ov S_{\De^R})^{2p-2}}
	\sqrt{\MO^a(\ov S_{\De^R})^{2}}.
	\eee
	Due to the convergence (\ref{to_S}) with $p=2$,
	we have $\MO^a(\ov S_{\De^R})^{2}\to 0$ as $R\to\infty$.
	Thus,
	it suffices to prove that the expectation
	$\MO^a |\ov S_{\De^R}|^{q}$
	is bounded uniformly in $R$, for any $q>0$.
	We have
	\bee\lbl{Sexp}
	|\ov S_{\De^R}|^q
	\leq C_q \big(
	e^{  \ov S_{\De^R}} +   e^{-\ov S_{\De^R}}
	\big).
	\eee
	Let us write
	\bee\non
	e^{\ov S_{\De^R}}
	= \wid \Psi_{\exp(\De^R)}.
	\eee
	Due to \rcor\ref{cor:mult}(2), we have
	\bee\lbl{cor2}
	\wid \Psi_{\exp(\De^R)} \ra \wid\Psi_{1} = 1 \ass
	R\ra\infty
	\qmb{in}\qu
	L^1(\CoR, \mP^a),
	\eee
	where $\wid\Psi_{1}$ is the multiplicative functional $\wid\Psi_{g}$
	corresponding to the function $g=1$.
	In particular, the $L^1$-norm
	$\MO^a e^{\ov S_{\De^R}}=\MO^a \wid \Psi_{\exp(\De^R)}$
	is bounded uniformly in $R$.
	Replacing $\Delta^R$ by $-\Delta^R$,
	the same argument implies that
	the expectation
	$\MO^a e^{-\ov S_{\De^R}}$ is also bounded uniformly in $R$.
	Then, due to (\ref{Sexp}),
	we see that the
	expectation $\MO^a |\ov S_{\De^R}|^q$
	is bounded uniformly in $R$ as well.
	So that, we obtain the desired convergence (\ref{to_S}).
	
	\ssk
	{\it Term $\ov\Psi_{b,a}^{R,\de} S_{h^\de_<}$.}
	Let us factorize
	\bee\lbl{Psi<}
	\ov\Psi_{b,a}^{R,\de}S_{h^\de_<}
	= \fr{\wid\Psi_{g_1^{R}}\wid\Psi_{g_2^\de}\Psi_{g_3^{\de}}S_{h^\de_<}}
	{\MO^a\big(\wid\Psi_{g_1^{R}}\wid\Psi_{g_2^\de}\Psi_{g_3^{\de}}\big)}
	\eee
	where
	$$
	g_1^{R}:= \Big( \Big(\fr{x-b}{x-a}\Big)^2-1 \Big) \mI_{\{x:\,|x|< R,\, |x-b|\geq 1\}} + 1
	$$
	and
	$$
	g_2^\de:= \Big( \fr{1}{(x-a)^2}-1 \Big) \mI_{\{x:\,|x-a|> \de,\, \de<|x-b|< 1\}} + 1,
	\quad
	g_3^\de:= \big( (x-b)^2-1 \big) \mI_{\{x:\,|x-a|> \de,\, \de<|x-b|< 1\}} + 1.
	$$
	Set
	$$
	g_1:= \Big( \Big(\fr{x-b}{x-a}\Big)^2-1 \Big) \mI_{\{x:\, |x-b|\geq 1\}} + 1
	$$
	and
	$$
	g_2:=\Big( \fr{1}{(x-a)^2}-1 \Big) \mI_{\{x:\,|x-b|< 1\}} + 1,
	\qu
	g_3:= \big( (x-b)^2-1 \big) \mI_{\{x:\, |x-b|< 1\}} + 1.
	$$
	Corollary \ref{cor:mult}(2) states that
	\bee\lbl{dg1}
	\wid\Psi_{g_1^R} \ra \wid\Psi_{g_1} \ass R\ra\infty
	\qu\qmb{in}\qu L^p(\CoR,\mP^a),
	\eee
	for any $p>0$, and
	\bee\lbl{dg2}
	\wid\Psi_{g_2^\de} \ra \wid\Psi_{g_2} \ass \de\ra 0 \qmb{in}\qu L^{1+\al}(\CoR,\mP^a),
	\eee
	for some $\al>0$.
	Since the functions $g_3^\de, g_3$ are bounded uniformly in $\de$ and $(g^\de_3-1),(g_3-1)$ have
	compact supports,
	we obviously have
	\bee\non
	\Psi_{g_3^\de} \ra \Psi_{g_3} \ass \de\ra 0
	\qmb{in}\qu L^p(\CoR,\mP^a),
	\eee
	for any $p>0$.
	Then the H\"older inequality implies
	\bee\non
	\MO^a \big(\wid\Psi_{g_1^R}\wid\Psi_{g_2^\de}\Psi_{g_3^\de}\big)
	\ra
	\MO^a \big(\wid\Psi_{g_1}\wid\Psi_{g_2}\Psi_{g_3}\big)
	\ass R\ra\infty,\qu \de\ra 0.
	\eee
	To prove that the nominator of (\ref{Psi<}) converges,
	we note that
	\bee
	\lbl{PsiS}
	\Psi_{g_3^\de} S_{h^\de_<}
	=2\sli_{\sck{x\in X: \\ |x-a|> \de,\, \de<|x-b|< 1}} (b-x)
	\prod\limits_{\sck{y\in X:\, y\neq x \\|y-a|>\de,\,\de<|y-b|< 1}}
	(y-b)^2.
	\eee
	Clearly, the right-hand side of (\ref{PsiS})
	converges as $\de\ra 0$ in $L^p(\CoR,\mP^a)$,
	for any $p>0$.
	Together with (\ref{dg1})-(\ref{dg2}),
	by the H\"older inequality this implies
	that the numerator of (\ref{Psi<}) converges in $L^1(\CoR,\mP^a)$
	as $R\ra\infty$, $\de\ra 0$,
	so that the function
	$\ov\Psi_{b,a}^{R,\de}S_{h_>^\de}$
	also does.
	\qed
	
	\bsk
	
	{\bf Acknowledgements.}
	A. Bufetov's research has received funding from the European Research Council
	(ERC) under the European Union's Horizon 2020 research and innovation
	programme (grant agreement No 647133 (ICHAOS)). It has also been
	funded by the Russian Academic Excellence Project `5-100',
	the grant MD 5991.2016.1 of the President of the Russian Federation
	and
	by the Gabriel Lam\'e Chair at the Chebyshev
	Laboratory of the SPbSU,
	a joint initiative of the French Embassy in the Russian Federation
	and the Saint-Petersburg State University.
	
	The work of A. Dymov was supported by the Russian Science Foundation under grant 14-21-00162 and performed in Steklov Mathematical Institute of RAS.
	In conformity with the reglementation of the Russian Science Foundation supporting the research of A. Dymov
	we  must indicate that A. Dymov prepared Subsections 3.2 and 3.3,  while A. Bufetov and H. Osada prepared Sections 1, 2, Subsection 3.1 and Appendix A.
	
	The research of H. Osada was supported by JSPS KAKENHI Grant Number JP16H06338.
	
	\appendix
	\numberwithin{equation}{section}

	\section{Regularization of additive and multiplicative functionals}
	\lbl{sec:regularization}
	
	In this appendix we
	consider a determinantal point process $\mP$ on $\mR$
	with the kernel $\Pi$ and assume that $\Pi$ satisfies
	Assumption \nolinebreak 2.
	We explain results from \cite{Buf}
	which we use in this paper
	and prove some auxiliary convergence results.
	
	\subsection{Additive functionals}
	\lbl{sec:add}
	
	Let $f:\mR \to {\mathbb C}$ be a Borel function.
	Define the corresponding {\it additive functional}
	\bee\lbl{addfunc}
	S_f:\CoR\mapsto\mR,
	\qu
	S_f(X)=\sum_{x\in X} f(x),
	\eee
	where the series may diverge.
	If
	$S_f\in L_1(\Conf(\mR), \mP)$, then we introduce the
	{\it normalized additive functional}
	\begin{equation}\label{normaddfunc}
	{\overline S}_f=S_f-\MO S_f.
	\end{equation}
	Now we will show that the normalized additive functional can be defined even when the additive functional itself is not well-defined.
	Introduce the Hilbert space $\VV(\Pi)$ of real functions
	with the norm
	\bee\non
	\|f\|_{\VV(\Pi)}^2=\displaystyle \frac12
	\ilif \ilif |f(x)-f(y)|^2 |\Pi(x,y)|^2dxdy.
	\eee
	Here we identify functions which differ by a constant.
	If a function $f$ is such that
	$S_f\in L^2(\CoR,\mP)$,
	we have
	$$
	\MO |\ov S_f|^2=\Var S_f=\|f\|_{\VV(\Pi)}^2.
	$$
	In particular, this is the case if the function $f$
	is bounded and has compact support.
	Thus, the
	correspondence $f\to {\overline S}_{f}$ is an isometric embedding of a dense subset of $\VV(\Pi)$ into $L_2(\Conf(\mR), \mP)$.
	It  therefore admits a unique  isometric extension onto the whole space $\HH$, and we get
	\bpp
	\lbl{prop:addfunc}
	There exists a unique linear isometric embedding
	$${\overline S}: \VV(\Pi) \hookrightarrow  L_2(\mathrm{Conf}(\mR), \mP),
	\qu
	{\overline S}: f\to {\overline S}_f
	$$
	such that
	\begin{enumerate}
		\item $\MO {\overline S}_f=0$ for all $f\in \VV(\Pi)$;
		\item if $S_f\in L_1(\mathrm{Conf}(\mR), \mP)$, then ${\overline S}_f$ is given by (\ref{normaddfunc}).
	\end{enumerate}
	\epp
	For more details see Proposition 4.1 in \cite{Buf}.
	
	\ssk
	Let $\mP^a$ be the reduced Palm measure of the measure $\mP$,
	conditioned at the point $a$.
	\begin{cor}\lbl{cor:add}
		\begin{enumerate}
			\item
			For $\hat{\rho}_1$-almost every $a\in\mR$ there exists a function $\ov S_a:\mathrm{Conf}(\mR)\mapsto\mR$
			such that the convergence (\ref{S1conv}) takes place,
			uniformly in $a\in\mR$ ranging in a compact set.
			\item
			The convergence (\ref{to_S}) takes place for $p=2$,
			uniformly in $a,b\in\mR$ ranging in a compact set.
		\end{enumerate}
	\end{cor}
	{\it Proof.}
	In \cite{ST1} it is proven that
	the reduced Palm measure $\mP^a$ coincides for
	$\hat{\rho}_1$-almost every $a\in\mR$
	with the determinantal measure associated
	with the kernel $\Pi^a$, given by
	\bee\lbl{Prkern}
	\Pi^a(x,y):=\Pi(x,y)-\fr{\Pi(x,a)\Pi(a,y)}{\Pi(a,a)} \qmb{if}\qu \Pi(a,a)\neq 0
	\eee
	and $\Pi^a(x,y):=\Pi(x,y)$ if $\Pi(a,a)=0$.
	It can be checked directly that the kernel $\Pi^a$ satisfies Assumption 2.
	Then item {\it 1} follows
	from \rprop \ref{prop:addfunc}
	applied to the kernel $\Pi^a$.
	Indeed, we have
	$$
	\ov S_a^{R,\de}=\ov S_{f_a^{R,\de}}
	\qu\mbox{with}\qu
	f_a^{R,\de}(x)=\fr{2}{a-x}\mI_{\{|x|<R,\,|x-a|>\de\}}.
	$$
	Clearly,
	$f_a^{R,\de}\ra f_a$ in $\VV(\Pi^a)$ uniformly in $a$ as $R\ra\infty$, $\de\ra 0$,
	where
	$$
	f_a(x):=\fr{2}{a-x}.
	$$
	Then
	\rprop \ref{prop:addfunc}
	implies the desired convergence with
	$\ov S_a:=\ov S_{f_a}$.
	Item {\it 2} can be proven by the same argument.
	\qed

	\subsection{Multiplicative functionals}
	\lbl{sec:mult}
	
	For a bounded nonnegative function $g$
	with compact support we define the
	{\it multiplicative functionals }
	$\Psi_g,\wid\Psi_g\,:\CoR\mapsto\mR$
	as
	\bee\lbl{multfunc}
	\Psi_g=\prod_{x\in X}g(x)=e^{S_{\log g}}
	\qnd
	\wid\Psi_g=e^{\ov S_{\log g}}.
	\eee
	Here we set $\Psi_g(X)=\wid \Psi_g(X)=0$
	if there is $x\in X$ such that $g(x)=0$.
	In view of \rprop\ref{prop:addfunc},
	we can extend the multiplicative functional
	$\wid\Psi_g$ to functions $g$
	satisfying $\|\log g\|_{\VV(\Pi)}<\infty$.
	If $\wid\Psi_g\in L^1(\CoR,\mP)$ we define the
	{\it normalized multiplicative functional} as
	\bee\lbl{normmultfunc}
	\ov\Psi_g=\fr{\wid\Psi_g}{\MO \wid\Psi_g}.
	\eee
	Fix positive  numbers $\alpha>0$, $\eps>0$, $M>\eps$ and two bounded Borel subsets $B^1,B^2\in\BB(\mR)$ satisfying
	$$
	|| \mI_{B^1 \cup B^2}\PPP ||<1.
	$$
	Denote by  $\mathscr G$ the set of nonnegative Borel functions 
	$g:\,\mR\mapsto\mR$ satisfying
	\begin{enumerate}
		\item $\{x: g(x)<\varepsilon\}\subset B^1$;
		\item $\{x: g(x)>M\}\subset B^2$;
		\item  ${\ds\int_{B^2} |g(x)|^{1+\alpha}\Pi(x,x)dx+ \int_{\mR\setminus B^2} |g(x)-1|^2\Pi(x,x) dx<\infty}$.
	\end{enumerate}
	We metrize $\mathscr G$ by equipping it with the distance
	$$
	d_{{\mathscr{G}}}(g_1, g_2)=\int_{B^2} |g_1(x)-g_2(x)|^{1+\alpha}\Pi(x,x)dx+ \int_{\mR\setminus B^2} |g_1(x)-g_2(x)|^2\Pi(x,x) dx.
	$$
	Then $\mathscr G$ becomes a complete separable metric space.
	Below we formulate Proposition \nolinebreak 4.3 from \cite{Buf}.
	\bpp\label{prop:mult}
	For any $\alpha^{\prime}: 0<\alpha^{\prime}<\alpha$, the correspondences $g\to \widetilde \Psi_g$,  $g\to{\overline \Psi}_g$
	induce continuous mappings from $\mathscr G$ to $L_{1+\alpha^{\prime}}(\mathrm{Conf}(\mR), \mP)$.
	\epp
	\begin{cor}\lbl{cor:mult}
		\begin{enumerate}
			\item
			Assertion of Lemma \ref{lem:multf}(1) is satisfied.
			\item
			For $\hat{\rho}_1$-almost all $a\in\mR$,
			convergences (\ref{cor2}), (\ref{dg1}) and (\ref{dg2}) take place and are uniform in $a,b\in\mR$ as $a$ range in a compact set and $b$ satisfies $|b-a|<\theta\ll 1$.
		\end{enumerate}
	\end{cor}
	{\it Proof.}
	{\it Item 1}. As we have already explained in the proof of \rcor\ref{cor:add},
	the reduced Palm measure $\mP^a$ coincides for $\hat{\rho}_1$-almost all $a\in\mR$
	with the determinantal measure associated with the kernel $\Pi^a$
	given by \eqref{Prkern}.
	Moreover, $\Pi^a$ satisfies Assumption 2.
	Indeed, one can check that $\Pi^a$ is an orthogonal projection kernel
	onto the space $L^a\subset L^2(\mR,dx)$ defined as 
	$L^a:=\{\ph\in L:\, \ph(a)=0\}$, see Section 2.16 of \cite{Buf}.
	Then Assumptions \nolinebreak 2(1,2) follow. Assumptions 2(3,4) are obvious.
	Let
	$$
	g_{a,b}(x):=\Big(\fr{x-b}{x-a}\Big)^2
	\qnd
	g_{a,b}^{R,\de}(x):=(g_{a,b}(x)-1)\mI_{\{|x|>R,\, |x-b|>\de,\,|x-a|>\de\}}+1.
	$$
	Using that the function $\Pi^a_{diag}(x):=\Pi^a(x,x)$ 
	has zero of second order at the point $x=a$,
	we find that
	$g_{a,b},g^{R,\de}_{a,b}\in\mathscr G$,
	for an appropriate choice of the sets $B^1,B^2$ and numbers $\al,M$ (independent from $R,\de$ and $a,b$),
	where the space $\mathscr G$ is defined with respect to the kernel $\Pi^a$.
	Moreover,
	\bee\lbl{Gconv}
	d_{\mathscr G}(g^{R,\de}_{a,b},g_{a,b})\ra 0 \ass R\ra\infty,\de\ra 0,
	\eee
	uniformly in $a,b$ as they range in compact sets.
	Then \rprop\ref{prop:mult} implies that
	$$
	\ov\Psi_{g^{R,\de}_{a,b}}\to \ov\Psi_{g_{a,b}} \ass R\ra\infty, \de\ra 0
	\qu\mbox{in}\quad L^{1+\al'}(\CoR,\mP^a),
	$$
	for any $0<\al'<\al$.
	Since $\ov\Psi_{a,b}^{R,\de}=\ov\Psi_{g_{a,b}^{R,\de}}$,
	we get the desired convergence with $\ov\Psi_{a,b}=\ov\Psi_{g_{a,b}}.$
	Its uniformity in $a,b$ follows from the uniformity of convergence
	(\ref{Gconv}) by a direct analysis of the proof of \rprop\ref{prop:mult}
	(i.e. of \rprop 4.3 from \cite{Buf}).
	
	{\it Item 2.} Similarly with the last item, the desired convergences follow by applying \rprop\ref{prop:mult}
	with the kernel $\Pi^a$.
	The convergence \eqref{dg1} takes place for arbitrary $p>0$ since
	we assume $|a-b|\ll 1$, so that
	the functions $g^{R}_1$ and $g_1$ are bounded,
	and then $\al$ can be chosen arbitrarily large.
	\qed

\end{document}